\documentclass[a4paper,10pt]{amsart}
\usepackage[utf8]{inputenc}
\usepackage{amsmath,amsthm,amssymb,bm}
\usepackage{hyperref,mathtools,dsfont}
\usepackage{a4wide}
\usepackage{aliascnt}

\newtheorem{theorem}{Theorem}[section]
\newtheorem*{theorem*}{Theorem}

\newaliascnt{lemma}{theorem}

\aliascntresetthe{lemma}

\newaliascnt{fact}{theorem}

\aliascntresetthe{fact}

\newaliascnt{corollary}{theorem}

\aliascntresetthe{corollary}

\newaliascnt{proposition}{theorem}
\newtheorem{proposition}[proposition]{Proposition}
\aliascntresetthe{proposition}

\newaliascnt{df}{theorem}
\theoremstyle{definition}\newtheorem{df}[df]{Definition}
\aliascntresetthe{df}

\newaliascnt{remark}{theorem}
\theoremstyle{remark}\newtheorem{remark}[remark]{Remark}
\aliascntresetthe{remark}

\numberwithin{equation}{section}

\DeclareMathOperator{\Pol}{\mathcal{O}}
\DeclareMathOperator{\id}{id}

\newcommand{\GG}{\mathbb{G}}
\newcommand{\HH}{\mathbb{H}}

\newcommand{\XX}{\mathbb{X}}

\title{Quantum Bernstein's Theorem and the Hyperoctahedral Quantum Group}
\author{Pawe{\l} J\'oziak}
\address{Institute of Mathematics of the Polish Academy of Sciences,
ul.~\'Sniadeckich 8, 00--656 Warszawa, Poland and  Institute of mathematics, University of Wroc{\l}aw, pl. Grunwaldzki 2/4, 50--384 Wroc{\l}aw, Poland}
\email{pjoziak@impan.pl}

\author{Kamil Szpojankowski} 
\address{Faculty of Mathematics and Information Science, Warsaw University of Technology, ul. Koszykowa 75, 00-662 Warszawa, Poland.}
\email{k.szpojankowski@mini.pw.edu.pl}



\begin{document}
\begin{abstract}
We study an extension of Bernstein's theorem to the setting of quantum groups. For a $d$-tuple of free, identically distributed random variables we consider a problem of preservation of freeness under the action of a quantum subset of the free orthogonal quantum group. For a subset not contained in the hyperoctahedral quantum group we prove that preservation of freeness characterizes Wigner's semicircle law. We show that freeness is always preserved if the quantum subset is contained in the hyperoctahedral quantum group. We provide examples of quantum subsets which show that our result is an extension of results known in the literature.
\end{abstract}
\maketitle
\section*{Introduction}
The well-known classical Bernstein's theorem states that given a $d$-dimensional random vector $\underline{X}=(X_1,\ldots,X_d)$ with independent coordinates and a generic orthogonal transformation  (with the exception of those which are equivalent to signed permutations) $A\in O_d$, if the coordinates of the random vector $A\underline{X}$ are still independent, then the random vector necessarily consists of identically distributed Gaussian random variables. 

This result forms one of the cornerstones of modern theory of random matrices. Indeed, if one asks for an ensemble of random matrices, with independent entries on the upper triangular part, to be invariant under a change of basis (unitary in the complex case and orthogonal in the real case), then one is forced to deal with Gaussian ensembles, due to Bernstein's theorem. This can be seen as a no-go result: there is no generalization of GOE to non-Gaussian orthogonal ensembles. 

A result similar in nature to Bernstein's theorem was obtained by Nica: he showed in \cite{Nica96} that if a vector consists of free random variables that remain free after application of a generic rotation, then elements of this vector form necessarily a semicircular system (free probability analogue of Gaussian random variables).

The aim of this note is to push this study further. We ask whether a similar result can be obtained if one replaces orthogonal transformations $O_d$ by free quantum orthogonal transformations $O_d^+$ of Wang and van Daele \cite{vDW}, and what is the right analogue of \emph{generic orthogonal transformations} in this context, i.e. which \emph{transformations} do not characterize Wigner's law. This result can be seen as a sequel of series of recent works on applications of quantum group symmetry in the theory of free probability, originating from K\"ostler and Speicher's free de Finetti theorem \cite{KS09} (see also a survey \cite{Spe17}) and later developed by Curran and others \cite{BBC11,BCS11,BCS12,Cur09,Cur10,Cur11,CS11}.

Let $\Pol(O_d^+)$ be the Hopf ${}^*$-algebra of a free orthogonal quantum group $O_d^+$ and let $\XX\subset O_d^+$ be its subset, i.e. a ${}^*$-epimorphism $\beta\colon \Pol(O_d^+)\to\mathcal{B}=\Pol(\XX)$ for some ${}^*$-algebra $\mathcal{B}$. Let $(\mathcal{A},\varphi)$ be a non-commutative ${}^*$-probability space and consider a family of random variables $X_1,\ldots,X_d\in\mathcal{A}$. We will study the \emph{quantum family of rotated random variables by the transformations from $\XX$}: the random variables $Y_j=\sum_{i=1}^dX_i\otimes \dot{u}_{ij}\in\mathcal{A}\otimes\Pol(\XX)$, where $u=(u_{ij})_{1\leq i,j\leq d}$ is the fundamental corepresentation of $\Pol(O_d^+)$ and $\dot{u}_{ij}=\beta(u_{ij})$. Let us also denote by $H_d^+$ the hyperoctahedral quantum group \cite{BBC07}. The main result of this note states the following.

\begin{theorem*}[Quantum Bernstein's theorem]
Assume that $X_1,\ldots,X_d$ are free and identically distributed. If $Y_1=\sum_{i=1}^dX_i\otimes \dot{u}_{i1},\ldots,Y_d=\sum_{i=1}^dX_i\otimes \dot{u}_{id}$ are free with amalgamation over $\Pol(\XX)$ and $\XX\not\subset H_d^+$, then $X_1,\ldots,X_d$ form a semicircular system. Conversely, given a $d$-tuple of free random variables $X_1,\ldots,X_d\in(\mathcal{A},\varphi)$, the $d$-tuple $Y_j=\sum_{i=1}^dX_i\otimes\dot{u}_{ij}\in\mathcal{A}\otimes\Pol(H_d^+)$ is free with amalgamation over $\Pol(H_d^+)$.
\end{theorem*}

Let us stress that, unlike the aforementioned articles on quantum groups as the sources of distributional symmetries in free probability, we do not assume invariance of the joint distribution. We only assume that the freeness property is preserved, the joint distribution may, a priori, change. The plan of the manuscript is as follows: we gather some preliminaries in \autoref{sec:preliminaries}: we recall relevant notions from the theory of quantum groups in \autoref{sec:quantumgroups} and from operator-valued free probability theory in \autoref{sec:freeproba}; in particular we recall the notion of operator-valued free cumulants, one of the main tools which we use in the proof of the main result. Here is the point where we were not able to adopt in a straightforward way the proof of non-commutative Bernstein's theorem - there is no notion of multidimensional $R$-transform in the operator-valued setting, for more details see \cite{SpeCombOpVal}. We are forced to use operator-valued free cumulants instead. \autoref{sec:qbt} is the core of the paper: we start with the description of the \emph{change of coordinates} in free cumulants under the action of the free orthogonal quantum group. Next we prove the first part of the main result stated above. The aim of \autoref{sec:other} is to discuss optimality (\autoref{sec:preservation}) and non-triviality (\autoref{sec:nontriviality}) of our quantum Bernstein's theorem. Namely, we show that a random vector consisting of free entries remains free after application of rotations from the hyperoctahedral quantum group, whatever the starting marginal distributions were. This proves second part of the Theorem above. If any quantum subset $\XX\subset O_d^+$ contained classical \emph{generic rotation}, our result would be reducible to Nica's non-commutative Bernstein's Theorem. In \autoref{sec:nontriviality} we provide examples of quantum subsets $\XX\subset O_d^+$ satisfying assumptions of our quantum Bernstein's Theorem, yet having not enough points for Nica's result to be applicable.

\section{Preliminaries}\label{sec:preliminaries}
\subsection{Free orthogonal and hyperoctahedral quantum groups}\label{sec:quantumgroups}
The theory of compact quantum groups in operator-algebraic setting was initiated by Woronowicz \cite{SLW87a} (see \cite{SLW95b,timmermann} for more details), we briefly introduce key concepts of this theory. A unital $C^{\ast}$-algebra $\mathnormal{A}$ endowed with a ${}^{\ast}$-homomorphism $\Delta\colon \mathnormal{A}\to \mathnormal{A}\otimes \mathnormal{A}$ (the minimal tensor product) satisfying the coassociativity condition: $(\Delta\otimes\id)\Delta=(\id\otimes\Delta)\Delta$ is called a \emph{Woronowicz algebra}, if the \emph{cancellation laws} holds:
\[(\mathds{1}\otimes \mathnormal{A})\cdot\Delta(\mathnormal{A})=\mathnormal{A}\otimes \mathnormal{A} = (\mathnormal{A}\otimes\mathds{1})\cdot\Delta(\mathnormal{A}).\]
The Woronowicz algebra $\mathnormal{A}$ can be always endowed with a unique state $h\in \mathnormal{A}^{\ast}$, called the \emph{Haar state}, which is left and right invariant:
\[(\id\otimes h)\Delta=(h\otimes\id)\Delta=h(\cdot)\mathds{1}.\]

Such an algebra correspond to a \emph{compact quantum group} $\mathbb{G}$ via abstract extension of the Gelfand-Naimark duality: $\mathnormal{A}=C(\mathbb{G})$, the algebra of continuous functions on $\mathbb{G}$. This algebra contains a unique dense Hopf ${}^{\ast}$-subalgebra $\Pol(\mathbb{G})$ (i.e. $\Delta\restriction_{\Pol(\GG)}\colon\Pol(\GG)\to\Pol(\GG)\otimes_{alg}\Pol(\GG)$); it is spanned by matrix coefficients of unitary representations of $\mathbb{G}$. This Hopf ${}^{\ast}$-algebra can have many different $C^{\ast}$-completions: the norm induced by GNS construction for $h$, the one coming from $A$ and the universal $C^*$-norm need not coincide. However, in our considerations only the Hopf-algebraic structure of $\Pol(\GG)$, and existence of a faithful Hilbert space representation of $\Pol(\GG)$, are relevant. 

We call $\mathbb{G}$ a compact matrix quantum group if $C(\mathbb{G})$ can be given a \emph{fundamental} corepresentation $u\in M_n(C(\mathbb{G}))=\mathnormal{B}(\mathbb{C}^n)\otimes C(\mathbb{G})$. Denoting $u_{ij}=(\langle e_i|\cdot|e_j\rangle\otimes\id)u$ for a fixed ONB $(e_i)_{1\leq i\leq n}\subset \mathbb{C}^n$ (with the standard inner product),  $u$ is called a fundamental corepresentation if: 
\[\Delta(u_{ij})=\sum_{k=1}^n u_{ik}\otimes u_{kj}\qquad i,j\in\{1,2\ldots n\}\]
and \[\langle \{u_{ij}:1\leq i,j\leq n\}\rangle=\Pol(\mathbb{G})\]
where $\langle X\rangle$ denotes the ${}^{\ast}$-algebra generated by elements of $X\subseteq C(\mathbb{G})$. Then the Hopf-algebraic structure on $\Pol(\GG)$ is uniquely determined by the prescription $S(u_{ij})=u_{ji}^*$ and $\varepsilon(u_{ij})=\delta_{ij}$. 

Sticking to the Gelfand-Najmark picture, a homomorphism $\HH\to\GG$ is, by definition, the transpose of a ${}^*$-homomorphism $\pi\colon \Pol(\GG)\to\Pol(\HH)$ which intertwine the coproducts: $(\pi\otimes\pi)\circ\Delta_{\GG}=\Delta_{\HH}\circ\pi$. If moreover $\pi$ is a surjection, then we say that $\HH$ is a subgroup of $\GG$ and write $\HH\subset\GG$. 

In this article we will be mainly interested in three examples of compact matrix quantum groups: the free orthogonal quantum group $O_d^+$, the hyperoctahedral quantum group $H_d^+$ and the quantum group of symmetries of a cube $O_{-1}(d)$.

\begin{df}[{\cite{vDW}}]
 Consider the universal $C^*$-algebra $C^u(O_d^+)$ generated by $d^2$ generators $u_{ij}$, $1\leq i,j\leq d$ subject to the relations:
 \begin{enumerate}
  \item all generators are self-adjoint $u_{ij}=u_{ij}^*$;
  \item the matrix $u=(u_{ij})_{1\leq i,j\leq d}$ is orthogonal, i.e. $u^{\top}u=uu^{\top}=\mathds{1}\in M_d(C^u(O_d^+))$.
 \end{enumerate}
 The $C^*$-algebra $C^u(O_d^+)$ is an algebra of continuous functions on a compact quantum group $O_d^+$, where the group-structure on $O_d^+$ is given the fundamental corepresentation $u$. The quantum group $O_d^+$ is called \emph{the free orthogonal quantum group}. 
 \end{df}
 For later reference, let us unpack condition (2):
 \begin{equation}\label{eq:normone}
  \sum_{i=1}^du_{ij}^2=\mathds{1}=\sum_{j=1}^du_{ij}^2
 \end{equation}
 and 
\begin{equation}\label{eq:ortho}
  \sum_{i=1}^du_{ij}u_{ij'}=0=\sum_{j=1}^du_{ij}u_{i'j}
 \end{equation}
 for all non-quantified indices $i\neq i'$ or $j\neq j'$.
\begin{df}[{\cite{BBC07}}]
 Consider the universal $C^*$-algebra $C^u(O_{-1}(d))$ generated by $d^2$ generators $u_{ij}$, $1\leq i,j\leq d$ subject to the relations:
 \begin{enumerate}
  \item all generators are self-adjoint $u_{ij}=u_{ij}^*$;
  \item the matrix $u=(u_{ij})_{1\leq i,j\leq d}$ is orthogonal, i.e. $u^{\top}u=uu^{\top}=\mathds{1}\in M_d(C^u(O_{-1}(d)))$;
  \item $u_{ij}u_{kl}=u_{kl}u_{ij}$ unless $i=k$ or $j=l$ (that is, unless these element lie in the same column or row in the matrix $u$, they commute);
  \item $u_{ij}u_{ij'}=-u_{ij'}u_{ij}$ if $j\neq j'$, and $u_{ij}u_{i'j}=-u_{i'j}u_{ij}$ if $i\neq i'$ (that is, different elements in the same row/column anticommute).
 \end{enumerate}
 The $C^*$-algebra $C^u(O_{-1}(d))$ is an algebra of continuous functions on a compact quantum group $O_{-1}(d)$, where the group-structure on $O_{-1}(d)$ is given by the fundamental corepresentation $u$. 
 \end{df}
\begin{df}[{\cite{BBC07}}]
 Consider the universal $C^*$-algebra $C^u(H_d^+)$ generated by $d^2$ generators $u_{ij}$, $1\leq i,j\leq d$ subject to the relations:
\begin{enumerate}
 \item all generators are self-adjoint $u_{ij}=u_{ij}^*$;
 \item the matrix $u=(u_{ij})_{1\leq i,j\leq d}$ is orthogonal, i.e. $u^{\top}u=uu^{\top}=\mathds{1}\in M_d(C^u(H_d^+))$;
 \item for all $i,j$ one has $u_{ij}=u_{ij}^3$ (or, equivalently, $u_{ij}^2=u_{ij}^4$, or, equivalently, $\sigma(u_{ij})\subset\{\pm1,0\}$);
 \item[(3')] $u_{ij}u_{ij'}=0=u_{ji}u_{j'i}$ for all $i,j\neq j'$.
\end{enumerate}
 The $C^*$-algebra $C^u(H_d^+)$ is an algebra of continuous functions on a compact quantum group $H_d^+$, where the group-structure on $H_d^+$ is given by the fundamental corepresentation $u$. The quantum group $H_d^+$ is called \emph{the hyperoctahedral quantum group}. 
 \end{df}
\begin{remark}
Note that relations $(3)$ and $(3^\prime)$ are equivalent and it is redundant to include both of them in the definition. On the other hand, in some computations one or the other will be more useful, we put them both for convinience of the reader. The equivalence of $(3)$ and $(3^\prime)$ was stated in \cite[Proposition 11.4(3)]{BBCC}, and the sketch of the passage between the conditions is as follows: computing $(\sum_{i}u_{ij}^2)^2$, using \eqref{eq:normone}, $(3^\prime)$ and positivity of $u_{ij}(\mathds{1}-u_{ij}^2)u_{ij}$ one arrives at (3). Conversly, sum of projections is a projection if and only if they are mutually orthogonal, hence \eqref{eq:normone} and (3) implies $u_{ij}^2u_{ij'}^2=0$. Then the $C^*$-identity, self-adjointness and the standard formula for spectral radius (which is equal to the norm for a self-adjoint element) yield the equality of the left-hand side of $(3^\prime)$, the right-hand side follows from applying antipode $S$.
\end{remark}
\begin{remark}\label{rmk:chainofsubgroups}
 Observe that $H_d^+, O_{-1}(d)\subset O_d^+$. Let us denote the ideals (in $\Pol(O_d^+)$) generated by the relations defining $H_d^+$ and $O_{-1}(d)$ by $I_{H_d^+}$ and $I_{O_{-1}(d)}$, respectively, observe that $I_{O_{-1}(d)}\subset I_{H_d^+}$.
  Indeed, $I_{H_d^+}=\ll\{u_{ij}u_{ij'}:i,j\neq j'\}\cup\{u_{ij}u_{i'j}:i\neq i',j\}\gg$, whereas $I_{O_{-1}(d)}=\ll \{u_{ij}u_{ij'}+u_{ij'}u_{ij}:i,j\neq j'\}\cup \{u_{ij}u_{i'j}+u_{i'j}u_{ij}:i\neq i',j\}\cup\{u_{ij}u_{i'j'}-u_{i'j'}u_{ij}:i\neq i',j\neq j'\}\gg$, where $\ll X\gg$ denotes the ideal generated by $X$. Thus the canonical map $\Pol(O_d^+)\ni u_{ij}\mapsto u_{ij}\in\Pol(H_d^+)$ factor through $\Pol(O_d^+)\ni u_{ij}\mapsto u_{ij}\in\Pol(O_{-1}(d))$. In other words, we have $H_d^+\subset O_{-1}(d)\subset O_d^+$ canonically.
\end{remark}

\subsection{Non-commutative probability}\label{sec:freeproba}
Here we recall basic definitions and facts concerning free probability and its operator valued extension. For more details we refer to \cite{NS2006,SpeCombOpVal,MS16}.

\begin{df}
	A non-commutative probability space consists of a pair $(\mathcal{A},\varphi)$ where $\mathcal{A}$ is a unital $*$-algebra and $\varphi:\mathcal{A}\to \mathbb{C}$ is linear functional such that $\varphi(\mathds{1})=1$ and $\varphi(aa^*)\geq 0$ for any $a\in\mathcal{A}$.
\end{df}

In this paper we will deal only with compactly supported measures. Such measures are uniquely determined by their moment sequences and we will identify the distribution of a random variable via moments.
\begin{df}
	For a self-adjoint random variable $a\in\mathcal{A}$ the distribution of $a$ is the unique probability measure such that
	\[\varphi\left(a^n\right)=\int t^n d\mu(t),\]
	for all $n=1,2,\ldots$.
\end{df}
In the non-commutative setting it is possible to define new notions of independence. The most prominent non-commutative independence is freeness defined by Voiculescu in \cite{Voi86}.
\begin{df}
		Consider a NCPS $(\mathcal{A},\varphi)$ and a family of unital subalgebras $\left(\mathcal{A}_i\right)_{i\in I}$ . The subalgebras $\left(\mathcal{A}_i\right)_{i\in I}$ are free if $\varphi(a_1\cdots a_n)=0$ whenever $a_i\in \mathcal{A}_{j_i}$, $j_1\neq j_2\neq \ldots \neq j_n$ and $\varphi(a_i)=0$ for all $i=1,\ldots,n$ and $n=1,2,\ldots$. Similarly, self-adjoint random variables $a,\,b\in\mathcal{A}$ are free (freely independent) when subalgebras generated by $(a,\,\mathds{1})$ and $(b,\,\mathds{1})$ are freely independent.
\end{df}
It turns out that for free random variables an analogue of Central Limit Theorem holds, that is if one takes a sequence $(a_n)_{n\geq 1}$ of identically distributed, free random variables with mean zero and variance one, then the distribution of the sequence 
\[\frac{a_1+\ldots+a_n}{\sqrt{n}}\]
tends to a universal limit, which is Wigner's semicircle law which has the density
\begin{align*}
\frac{1}{2\pi}\sqrt{4-x^2}\mathds{1}_{[-2,2]}(x).
\end{align*}
Free random variables can be succinctly studied in terms of the so called free cumulants. 

Let $\chi=\{B_1,B_2,\ldots\}$ be a  partition of the set of numbers $\{1,\ldots,k\}$. A partition $\chi$ is a crossing partition if there exist distinct blocks $B_r,\,B_s\in\chi$ and numbers $i_1,i_2\in B_r$, $j_1,j_2\in B_s$ such that $i_1<j_1<i_2<j_2$. Otherwise $\chi$ is called a non-crossing partition. The set of all non-crossing partitions of $\{1,\ldots,k\}$ is denoted by $NC(k)$.

\begin{df}
For any $k=1,2,\ldots$, (joint) cumulants of order $k$ of non-commutative random variables $a_1,\ldots,a_n$ are defined recursively as $k$-linear maps $\kappa_k:\mathcal{A}^k\to\mathbb{C}$ through equation
$$
\varphi(a_1\cdot\ldots\cdot a_m)=\sum_{\pi\in NC(m)}\,\prod_{B\in\pi}\,\kappa_{|B|}(a_i,\,i\in B)
$$
with $|B|$ denoting the size of the block $B$. One can write the above as
\begin{align*}
\varphi(a_1\cdots a_n)=\sum_{\pi\in NC(n)} \kappa_{\pi}(a_1,\ldots,a_n)
\end{align*}

Cumulants of single variable $a$, are defined in the same manner as above, one takes $a_1=\ldots=a_n=a$. We denote $\kappa_n(a)=\kappa_n(a,\ldots,a)$.
\end{df}
With free cumulants in hand one can easily characterize the semicircle distribution; for details we refer to \cite[Lecture 11]{NS2006}.

\begin{proposition}
	A random variable $a$ has the semicircular distribution with mean $\mu$ and variance $\sigma^2$ if and only if $\kappa_1(a)=\mu,\kappa_2(a)=\sigma^2$ and $\kappa_n(a)=0$ for all $n\geq 3$.
\end{proposition}

Next we discuss briefly extension of the concept of non-commutative probability space - so called operator valued non-commutative probability space. Observe that when one takes $\mathcal{B}=\mathbb{C}\mathds{1}$, then we are in the scalar valued framework, discussed above.

\begin{df}
	A $\mathcal{B}$-valued non-commutative probability space (NCPS) consists of a triple $(\mathcal{A},\mathcal{B},E)$, where $\mathcal{B}\subset\mathcal{A}$ and $E$ is a conditional expectation.
\end{df}

\begin{df}
	For a ${}^*$-algebra $\mathcal{A}$ and its ${}^*$-subalgebra $\mathcal{B}$ a linear map $E:\mathcal{A}\to\mathcal{B}$ is called a conditional expectation if
	\begin{align*}
	&E(b)=b \qquad\forall b\in\mathcal{B}\\
	&\mbox{and} \\
	&E(b_1 a b_2)=b_1 E(a) b_2 \qquad\forall a\in \mathcal{A}\, \mbox{and} \,\forall b_1,b_2\in \mathcal{B}
	\end{align*}
\end{df}

\begin{df}
	Let $x_1,\ldots,x_d\in \mathcal{A}$, then the joint distribution of $x_1,\ldots,x_d$ is given by all joint moments of the form
	\[
	E(y_1b_1\ldots b_{n-1}y_n),
	\]
	where $y_i\in\{x_1,\ldots,x_d\}$ and $b_i\in\mathcal{B}$, $i=1,\ldots,n$ and $n\geq 1$.
\end{df}

\begin{df}
	Consider a $\mathcal{B}$-valued NCPS $(\mathcal{A},\mathcal{B},E)$ and a family of subalgebras $\left(A_i\right)_{i\in I}$ where $\mathcal{B}\subset A_i$ for all $i\in I$ . The subalgebras $\left(A_i\right)_{i\in I}$ are free with amalgamation over $\mathcal{B}$ if $E(a_1\cdots a_n)=0$ whenever $a_i\in A_{j_i}$, $j_1\neq j_2\neq \ldots \neq j_n$ and $E(a_i)=0$ for all $i=1,\ldots,n$ and $n=1,2,\ldots$.
\end{df}

\begin{df}
	For an operator valued NCPS we define the corresponding $\mathcal{B}$-valued cumulants $(\kappa^{\mathcal{B}}_n)_{n\geq 1}$ via the moment-cumulant formula
	\begin{align*}
	E(a_1\cdots a_n)=\sum_{\pi\in NC(n)} \kappa_{\pi}^{\mathcal{B}}(a_1,\ldots,a_n)
	\end{align*}
	where cumulants are nested inside each other according to the nesting of blocks of $\pi$.
\end{df}
\begin{remark}
	For our purposes it is important that one can write an explicit formula for cumulants in terms of moments
	\begin{align} \label{MC_formula}
		\kappa_{n}^{\mathcal{B}}(a_1,\ldots,a_n)=\sum_{\pi\in NC(n)} E_{\pi}(a_1,\ldots,a_n) \mu(\pi,\mathit{1}_n),
	\end{align}
	where $E_\pi$ is a multiplicative functional defined on the lattice $NC(n)$ and again moments are  nested according to the nesting of blocks of $\pi$ (for the explanation of nesting, see \cite[p. 240]{MS16}). By $\mu(\cdot,\cdot)$ we mean the M\"obius function on $NC(n)$ and $\mathit{1}_n$ is the maximal partition in $NC(n)$ with respect to the reversed refinement order.

	\end{remark}

The relevance of operator valued cumulants stems from the following result proved in \cite{SpeCombOpVal}.

\begin{theorem} \label{thm:freeCumVanish}
	Let $(\mathcal{A},\mathcal{B}, E)$ be a $\mathcal{B}$-valued probability space and $(x_i)_{i\in I}$ a
	family of random variables in $\mathcal{A}$. Then the family $(x_i)_{i\in I}$ is free with amalgamation over $\mathcal{B}$ if and only if
	\begin{align*}
	\kappa_n^{\mathcal{B}}(y_1 b_1,\ldots,y_n b_n)=0
	\end{align*}
	for every $n\geq2$, all $b_1,\ldots,b_n\in \mathcal{B}$ and every non-constant choice of $y_1,\ldots,y_n\in\{x_i:i\in I\}$
\end{theorem}

\section{Quantum orthogonal transformations applied to random vectors}\label{sec:qbt}
Assume $\XX\subset O_d^+$ is a closed quantum subset of quantum orthogonal group, i.e. there is a surjective ${}^*$-homomorphism $\beta\colon \Pol(O_d^+)\twoheadrightarrow\Pol(\XX)$ (in particular, $\Pol(\XX)$ is generated, as an algebra, by elements $\beta(u_{ij})$, which we denote later by $\dot{u}_{ij}$). 

Let $(\mathcal{A},\varphi)$ be an NCPS. Then on $\mathcal{A}\otimes \Pol(\XX)$ one can define a conditional expectation given by $E=\varphi \otimes\id$, and we obtain an operator valued NCPS $(\mathcal{A}\otimes \Pol(\XX), \Pol(\XX),E)$.

\begin{remark} (Moments and cumulants of rotated variables)
	Take random variables $X_1,\ldots,X_d\in\mathcal{A}$, and define $Y_j=\sum_{i=1}^dX_i\otimes \dot{u}_{ij} \in \mathcal{A}\otimes \Pol(\XX)$.
	Then by the definition of $E$ and the fact that $\mathcal{A}$ and $\Pol(\XX)$, seen as subalgebras in $\mathcal{A}\otimes\Pol(\XX)$, commute, we have
	\begin{align*}
	E(Y_{j_1}\cdots Y_{j_n})=\sum_{i_1,\ldots,i_n=1}^d \varphi(X_{i_1},\ldots,X_{i_n})\dot{u}_{i_1 j_1}\cdots \dot{u}_{i_n j_n}.
	\end{align*}
	From the above we get that 
	\begin{align*}
	E_{\pi}(Y_{j_1}\cdots Y_{j_n})=\sum_{i_1,\ldots,i_n=1}^d \varphi_{\pi}(X_{i_1},\ldots,X_{i_n})\dot{u}_{i_1 j_1}\cdots \dot{u}_{i_n j_n},
	\end{align*}
	which in turn, by formula \eqref{MC_formula} (and simple change of order of summation) implies 
	\begin{align} \label{cum_trans}
	\kappa_{n}^{\Pol(\XX)}(Y_{j_1},\ldots, Y_{j_n})=\sum_{i_1,\ldots,i_n=1}^d\kappa_{n}(X_{i_1},\ldots,X_{i_n})\dot{u}_{i_1 j_1}\cdots \dot{u}_{i_n j_n},
	\end{align}
	where $\kappa_{n}^{\Pol(\XX)}$ denote the operator-valued cumulants corresponding to the conditional expectation $E$ and $\kappa_{n}$ are scalar valued cumulants corresponding to $\varphi$.
\end{remark}

\begin{theorem}
	Assume that $X_1,X_2,\ldots,X_d$ are free and $Y_1=\sum_{i=1}^dX_i\otimes \dot{u}_{i1},Y_2=\sum_{i=1}^dX_i\otimes \dot{u}_{i2},\ldots,Y_d=\sum_{i=1}^dX_i\otimes \dot{u}_{id}$ are free with amalgamation over $\Pol(\XX)$. Assume moreover $X_i$ are identically distributed. Then $X_1, X_2,\ldots, X_d$ are semicircular random variables, unless $\XX\subset H_d^+$ (i.e. the homomorphism $\beta$ factors as $\Pol(O_d^+)\to\Pol(H_d^+)\to\Pol(\XX)$).
\end{theorem}
\begin{remark}
 We expect that the assumption of identical distribution of $X_i$ could be dropped, as in the classical Bernstein's theorem. However, with the proof we found, it is essential.
\end{remark}

\begin{proof}
The strategy of the proof is the following: freeness of $Y_1,Y_2,\ldots,Y_d$ translates into relations about $\dot{u}_{ij}$'s with coefficients $\kappa_n(X_i)$ by means of \autoref{thm:freeCumVanish}. We show that, if these coefficients do not vanish for $n\geq 3$, then $\dot{u}_{ij}^2=\dot{u}_{ij}^4$, the defining relation for the hyperoctahedral quantum group, which violates the assumption $\XX\not\subset H_d^+$.

We start with even cumulants: assume $n\geq4$ is even. For $j\neq j'$, by freeness of $Y_j$ and $Y_{j'}$, \autoref{thm:freeCumVanish} yields: 
\[\kappa_n^{\Pol(\XX)}(Y_{j'},Y_{j'},Y_j,\ldots,Y_j)=0\]
where $Y_{j'}$ appears only in first two spots. By \eqref{cum_trans}, this is equivalent to
\[\sum_{i_1,\ldots,i_n=1}^d\kappa_{n}(X_{i_1},\ldots,X_{i_n})\dot{u}_{i_1 j'} \dot{u}_{i_2 j'}\dot{u}_{i_3 j}\cdots \dot{u}_{i_n j}=0,\]
which simplifies, thanks to freeness of $X_1, X_2,\ldots,X_d$, to the formula
\[\sum_{i=1}^d\kappa_{n}(X_{i})\dot{u}_{i j'}^{2} \dot{u}_{i j}^{n-2}=0.\]
This formula is valid for any pair of indices $j\neq j'$. Sum them all (over $j'\neq j$) to obtain (with the aid of \eqref{eq:normone}):
\begin{equation}\label{eq:even}
\sum_{i=1}^d\kappa_{n}(X_{i})(\sum_{j'\neq j}\dot{u}_{i j'}^2)\dot{u}_{i j}^{n-2} =\sum_{i=1}^d\kappa_{n}(X_{i})(\mathds{1}-\dot{u}_{i j}^2)\dot{u}_{i j}^{n-2}=0. 
\end{equation}

Now as $\kappa_n(X_i)=\kappa_n(X_1)$, either $\kappa_n(X_i)=0$ or, by dividing \eqref{eq:even} by $\kappa_n(X_1)$ and rearranging the terms, we have that:
\[\sum_{i=1}^d\dot{u}_{ij}^{n/2-1}(\mathds{1}-\dot{u}_{i j}^2)\dot{u}_{i j}^{n/2-1}=0\]

The above formula is valid for any $j$. But as it is a sum of positive operators, this can happen only if each of the summands is a zero operator, thus for all $i,j$ we have that
\[\dot{u}_{ij}^n=\dot{u}_{ij}^{n+2}\]
which, by spectral calculus, implies that $\sigma(\dot{u_{i,j}})\subset\{\pm1,0\}$, and consequently $\dot{u}_{ij}^2=\dot{u}_{ij}^4$, as desired.

We now proceed with the case $n\geq3$ odd. The beginning follows similar line of argument. For $j\neq j'$, by freeness of $Y_j$ and $Y_{j'}$, \autoref{thm:freeCumVanish} yields: 
\[\kappa_n^{\Pol(\XX)}(Y_j,Y_{j'},Y_{j'},Y_{j_1},Y_{j_1},\ldots,Y_{j_{\frac{n-3}{2}}})=0\]
where $Y_j$ appears on the first spot, $Y_j$ appears only on the second and third spots, and each $Y_{j_1},\ldots, Y_{j_{\frac{n-3}{2}}}$ appears exactly twice on neighboring spots. By \eqref{cum_trans}, this is equivalent to
\[\sum_{i_1,\ldots,i_n=1}^d\kappa_{n}(X_{i_1},\ldots,X_{i_n})\dot{u}_{i_1 j} \dot{u}_{i_2 j'}\dot{u}_{i_3 j'}\cdots \dot{u}_{i_n j_{\frac{n-3}{2}}}=0,\]
which simplifies, thanks to freeness of $X_1, X_2,\ldots,X_d$, to
\[\sum_{i=1}^d\kappa_{n}(X_{i})\dot{u}_{i j} \dot{u}_{i j'}^2\dot{u}_{i j_1}^2\cdots \dot{u}_{i j_{\frac{n-3}{2}}}^2=0.\]
The above relation is valid for any pair of indices $j\neq j'$ and any choice of $j_1,\ldots, j_{\frac{n-3}{2}}$. After summing over all $j_1$, and then over all $j_2$, ..., and then over all $j_{\frac{n-3}{2}}$, and using \eqref{eq:normone}, we obtain:
\begin{equation}\label{eq:odd1}
\sum_{i=1}^d\kappa_{n}(X_{i})\dot{u}_{i j}\dot{u}_{i j'}^2 =0. 
\end{equation}
The above is valid for all pairs $j\neq j'$. We can now use this relation in two different ways. To simplify things, we use the assumption $\kappa_n(X_i)=\kappa_n(X_1)\neq0$ and divide the relations by this scalar. Firstly, summing over all $j'\neq j$, one obtains (similarly as in the case of $n$ even):
\[\sum_{i=1}^d\dot{u}_{i j}(\mathds{1}-\dot{u}_{i j}^2) =0,\]
or, equivalently, 
\begin{equation}\label{eq:odd2}\sum_{i=1}^d\dot{u}_{i j} =\sum_{i=1}^d\dot{u}_{i j}^3. \end{equation}
on the other hand, computing $X^*X$ where $X$ denotes the left hand side of \eqref{eq:odd1}, we arrive at:
\[\sum_{i,i'}\dot{u}_{ij'}^2\dot{u}_{ij}\dot{u}_{i'j}\dot{u}_{i'j'}^2=0\]
this relation is valid for all pairs $j\neq j'$. We divide the above sum into sum over $i=i'$ and over all pairs $(i,i')$ such that $i\neq i'$. Then sum over all $j\neq j'$ (with $j'$ fixed) and change the order of summation to obtain:
\[\sum_{i}\dot{u}_{ij'}^2(\sum_{j\neq j'}\dot{u}_{ij}^2)\dot{u}_{ij'}^2+\sum_{i,i':i\neq i'}\dot{u}_{ij'}^2(\sum_{j\neq j'}\dot{u}_{ij}\dot{u}_{i'j})\dot{u}_{i'j'}^2=0\]
where the middle sum of the first term is nothing but $\mathds{1}-\dot{u}_{ij'}^2$ thanks to \eqref{eq:normone}, and the sum in the middle of the second term is nothing but $-\dot{u}_{ij'}\dot{u}_{i'j'}$ thanks to \eqref{eq:ortho}. This can be rewritten as:
\[\sum_{i}\dot{u}_{ij'}^4=\sum_{i,i'}\dot{u}_{ij'}^3\dot{u}_{i'j'}^3=\left(\sum_{i}\dot{u}_{ij'}^3\right)^2\]
where we moved the order-six terms from the first sum to the right hand side. Now we use \eqref{eq:odd2} to change the right hand side. This yields:  
\[\sum_{i}\dot{u}_{ij'}^4=(\sum_{i}\dot{u}_{ij'})^2=\sum_{i}\dot{u}_{ij'}^2+\sum_{i,i':i\neq i'}\dot{u}_{ij'}\dot{u}_{i'j'}. \]
Now summing over all indices $j'$ yields:
\[\sum_{i,j'}\dot{u}_{ij'}^4=\sum_{i,j'}\dot{u}_{ij'}^2+\sum_{i,i':i\neq i'}\sum_{j'}\dot{u}_{ij'}\dot{u}_{i'j'}, \]
where the last term on the right is equal to $0$ thanks to \eqref{eq:ortho}. Thus
\[\sum_{i,j'}\dot{u}_{ij'}(\mathds{1}-\dot{u}_{ij'}^2)\dot{u}_{ij'}=0\]
is a combination of positive operators, which can be equal to zero only if each of them is equal to zero, which amounts to:
\[\dot{u}_{ij'}^2=\dot{u}_{ij'}^4\]
for all $i,j'$, as desired.
\end{proof}

\begin{remark}
 In the case $d=2$, one can avoid the assumption of identical distribution of $X_i$, by appropriately modifying the proof of Nica in the scalar case \cite[Theorem 5.1]{Nica96}. Indeed, the relation coming from $\kappa_n^{\Pol(\XX)}(Y_j,Y_j,Y_{j'},\ldots,Y_{j'})=0$ can be written as \[\kappa_n(X_1)u_{1j}^2u_{1j'}^{n-2}+\kappa_n(X_2)u_{2j}^2u_{2j'}^{n-2}=0.\]
 Use \eqref{eq:normone} and multiply both sides of the above equation by $u_{2j'}$ from the left to get:
 \[\kappa_n(X_1)u_{2j'}u_{1j'}^{n-2}+\kappa_n(X_2)u_{2j'}^{n-1}=\kappa_n(X_1)u_{2j'}u_{1j'}^{n}+\kappa_n(X_2)u_{2j'}^{n+1}.\]
 Thanks to \eqref{eq:ortho} this is equivalent to:
 \[\kappa_n(X_1)u_{2j'}u_{1j'}^{n-2}=\kappa_n(X_2)u_{2j'}^{n-1}(u_{2j'}^2-\mathds{1})-\kappa_n(X_1)u_{2j}u_{1j}u_{1j'}^{n-1}.\]
 Now $\kappa_n^{\Pol(\XX)}(Y_j,Y_{j'},\ldots,Y_{j'})=0$ translates into $\kappa_n(X_1)u_{1j}u_{1j'}^{n-1}=-\kappa_n(X_2)u_{2j}u_{2j'}^{n-1}$, whose left-hand side appears as the last term in the above displayed equation. Substituting it yields
 \[\kappa_n(X_1)u_{2j'}u_{1j'}^{n-2}=\kappa_n(X_2)u_{2j'}^{n-1}(u_{2j'}^2+u_{2j}^2-\mathds{1}).\]
 And the term in paranthesis on the right-hand side vanish due to \eqref{eq:normone}. Multiplying the equation by $u_{2j'}$ from the left and using \eqref{eq:normone} to $u_{2j'}^2$ on the left-hand side translates the relation to:
 \[\kappa_n(X_1)u_{1j'}^{n-2}=\kappa_n(X_1)u_{1j'}^n\]
 where now only a cumulant of a single variable appears. Similarly one gets remaining relations which involve only $\kappa_n(X_2)$ or with $u_{2j'}$. All of them lead to the conclusion that either squares of the elements are all projections, or all cumulants $\kappa_n(X_i)$ with $n\geq3$ vanish.
 
\end{remark}
\section{Relevance of quantum Bernstein's theorem}\label{sec:other}
 In classical Bernstein's theorem and in its free version one cannot take a rotation which only permutes and changes signs of coordinates of the vector, as then independence or freeness of coordinates is trivially preserved. In other words the hyperoctahedral group always preserve independence and freeness. In the quantum version the forbidden transformations come from hyperoctahedral quantum group. We show that this assumption was indeed necessary: the hyperoctahedral quantum group always preserves freeness of coordinates. On the other hand our result is closely related to Nica's free Bernstein's theorem, it is natural to ask whether it is a genuine extension. We provide here examples of quantum subsets for which Nica's theorem cannot be applied.\\

\subsection{Preservation of freeness under transformations from the hyperoctahedral quantum group}\label{sec:preservation}

Let $\underline{X}=(X_1,\ldots, X_d)$ be a tuple consistsing of pairwise free entries belonging to a single NCPS $(\mathcal{A},\varphi)$. Our aim is to show that, whatever the marginal distributions of $\underline{X}$ are, the tuple $\underline{Y}=(Y_1,\ldots,Y_d)$ obtained by $Y_j=\sum_{i=1}^dX_i\otimes\dot{u}_{ij}\in\mathcal{A}\otimes\Pol(H_d^+)$, is free with amalgamation over $\Pol(H_d^+)$ with respect to conditional expectation $E=\varphi\otimes\id$. To this end, we need to show that all mixed $\Pol(H_d^+)$-valued free cumulants vanish. We calculate
\[\kappa_n^{\Pol(H_d^+)}(Y_{j_1},\ldots,Y_{j_n})=\sum_{i_1,\ldots,i_n=1}^d\kappa_{n}(X_{i_1},\ldots,X_{i_n})\dot{u}_{i_1 j_1}\cdots \dot{u}_{i_n j_n}=\sum_{i=1}^d\kappa_{n}(X_i)\dot{u}_{ij_1}\cdots \dot{u}_{i j_n}=0\]
where we used, respectively, \eqref{cum_trans}, freeness of entries of $\underline{X}$, and condition $(3^\prime)$ from the definition of $H_d^+$ together with the fact that the sequence $j_1,\ldots,j_n$ is non-constant.
\subsection{Non-triviality of quantum Bernstein's theorem}\label{sec:nontriviality}
To show that our quantum Bernstein's theorem cannot be deduced from Nica's free Bernstein's theorem, we need to provide an example of a quantum subset $\mathbb{X}\subset O_d^+$ such that
\begin{enumerate}
 \item $\mathbb{X}\not\subset H_d^+$ (i.e. $\Pol(O_d^+)\to\Pol(\XX)$ cannot be factored through $\Pol(H_d^+)$);
 \item $\mathbb{X}\cap O_d^+\subset H_d$, i.e. the quotient of $\Pol(\XX)$ by the commutator ideal is a quotient of $\Pol(H_d)$. In other words, the only points in $\XX$ are those already in $H_d$.
\end{enumerate}
The first item shows that such an $\mathbb{X}$ satisfies the assumption of our quantum Bernstein's theorem, whereas (2) shows that there are too few points in $\XX$ for Nica's result to be applicable.

Let $\GG$ is a subgroup generated by $\XX$ (in the sense of \cite{BB10,SS16}). We expect that if $\XX$ satisfies assumptions of our quantum Bernstein's theorem (i.e. an identically distributed $d$-tuple of free random variables, after applying transformations from $\XX$, is free with amalgamation over $\Pol(\XX)$), then $\GG$ also satisfies these assumptions. However, we were not able to find a proof of this assertion without using the quantum Bernstein's theorem itself: under these assumptions, this tuple is necessarily a semicircular system, and Curran showed in \cite[Proposition 3.5]{Cur10} that its joint distribution remains unchanged after applying transformations from $O_d^+$, so in particular it is preserved by transformations from $\GG\subset O_d^+$. 

With $\mathbb{X}$ being a quantum subset, the theorem is formally stronger than if we assumed that $\XX$ is a subgroup. Note that the procedure of generation of quantum group from $\XX$ can yield points that are not obtained by generation of group from points of $\XX$, as discovered in \cite[Section 2.6]{pjrmk}. Because of that, and because of the statement from previous paragraph, it is desirable to provide an example of $\XX$ that is a quantum group from the very beginning.

The quantum group which we need is $O_{-1}(d)$. As noted in \autoref{rmk:chainofsubgroups}, we have that $H_d^+\subset O_{-1}(d)$, and by the results of \cite{BBC07}, we have that $H_d^+\neq O_{-1}(d)$ for $d\geq3$ (see \cite[Sections 5 \& 7]{BBC07}). It is straightforward to verify that the maximal classical subgroup of $O_{-1}(d)$ is precisely $H_d$: the conjunction of anticommutativity and commutativity in a single row/column forces entries to satisfy the relation $u_{ij}u_{ij'}=0=u_{ij}u_{i'j}$, a defining relation for $C(H_d)$.

Thus, to get $\XX$ as described in first paragraph of the subsection, it is enough to take (for $d\geq3$) any intermediate quantum space $H_d^+\subsetneq\mathbb{X}\subseteq O_{-1}(d)$, i.e. any intermediate quotient $\Pol(O_{-1}(d))\to\Pol(\XX)\to\Pol(H_d^+)$, with the latter arrow being a proper surjection. In particular one can take $\XX=O_{-1}(d)$, which is itself a quantum group.

\section*{Acknowledgement} The authors would like to thank Marek Bo{\.z}ejko and Jacek Weso{\l}owski for encouragement during development of this project. We are also immensely grateful to Adam Skalski for his comments on an earlier version of the manuscript. PJ was partially supported by the NCN (National Science Center) grant 2015/17/B/ST1/00085. KSz was partially supported by the NCN (National Science Center) grant 2016/21/B/ST1/00005. 

\bibliographystyle{alpha}
\bibliography{QBT}

\begin{thebibliography}{BBCC11}

\bibitem[BB10]{BB10}
Teodor Banica and Julien Bichon.
\newblock Hopf images and inner faithful representations.
\newblock {\em Glasg. Math. J.}, 52(3):677--703, 2010.

\bibitem[BBC07]{BBC07}
Teodor Banica, Julien Bichon, and Beno{\^\i}t Collins.
\newblock The hyperoctahedral quantum group.
\newblock {\em J. Ramanujan Math. Soc.}, 22(4):345--384, 2007.

\bibitem[BBC11]{BBC11}
Teodor Banica, Julien Bichon, and Stephen Curran.
\newblock Quantum automorphisms of twisted group algebras and free
  hypergeometric laws.
\newblock {\em Proc. Amer. Math. Soc.}, 139(11):3961--3971, 2011.

\bibitem[BBCC11]{BBCC}
T.~Banica, S.~T. Belinschi, M.~Capitaine, and B.~Collins.
\newblock Free {B}essel laws.
\newblock {\em Canad. J. Math.}, 63(1):3--37, 2011.

\bibitem[BCS11]{BCS11}
Teodor Banica, Stephen Curran, and Roland Speicher.
\newblock Stochastic aspects of easy quantum groups.
\newblock {\em Probab. Theory Related Fields}, 149(3-4):435--462, 2011.

\bibitem[BCS12]{BCS12}
Teodor Banica, Stephen Curran, and Roland Speicher.
\newblock De {F}inetti theorems for easy quantum groups.
\newblock {\em Ann. Probab.}, 40(1):401--435, 2012.

\bibitem[CS11]{CS11}
Stephen Curran and Roland Speicher.
\newblock Quantum invariant families of matrices in free probability.
\newblock {\em J. Funct. Anal.}, 261(4):897--933, 2011.

\bibitem[Cur09]{Cur09}
Stephen Curran.
\newblock Quantum exchangeable sequences of algebras.
\newblock {\em Indiana Univ. Math. J.}, 58(3):1097--1125, 2009.

\bibitem[Cur10]{Cur10}
Stephen Curran.
\newblock Quantum rotatability.
\newblock {\em Trans. Amer. Math. Soc.}, 362(9):4831--4851, 2010.

\bibitem[Cur11]{Cur11}
Stephen Curran.
\newblock A characterization of freeness by invariance under quantum spreading.
\newblock {\em J. Reine Angew. Math.}, 659:43--65, 2011.

\bibitem[J{\'o}z]{pjrmk}
Pawe{\l} J{\'o}ziak.
\newblock Remarks on hopf images and quantum permutation groups ${S}_n^+$.
\newblock {\em Canad. Math. Bull.}
\newblock In press.

\bibitem[KS09]{KS09}
Claus K\"ostler and Roland Speicher.
\newblock A noncommutative de {F}inetti theorem: invariance under quantum
  permutations is equivalent to freeness with amalgamation.
\newblock {\em Comm. Math. Phys.}, 291(2):473--490, 2009.

\bibitem[MS17]{MS16}
James~A. Mingo and Roland Speicher.
\newblock {\em Free Probability and Random Matrices}.
\newblock Springer, 2017.

\bibitem[Nic96]{Nica96}
Alexandru Nica.
\newblock {$R$}-transforms of free joint distributions and non-crossing
  partitions.
\newblock {\em J. Funct. Anal.}, 135(2):271--296, 1996.

\bibitem[NS06]{NS2006}
Alexandru Nica and Roland Speicher.
\newblock {\em Lectures on the combinatorics of free probability}, volume 335
  of {\em London Mathematical Society Lecture Note Series}.
\newblock Cambridge University Press, Cambridge, 2006.

\bibitem[Spe98]{SpeCombOpVal}
Roland Speicher.
\newblock Combinatorial theory of the free product with amalgamation and
  operator-valued free probability theory.
\newblock {\em Mem. Amer. Math. Soc.}, 132(627):x+88, 1998.

\bibitem[Spe14]{Spe17}
Roland Speicher.
\newblock Free probability and non-commutative symmetries.
\newblock In {\em Quantum Symmetries (Metabief)}. Springer, 2014.
\newblock In preparation.

\bibitem[SS16]{SS16}
Adam {Skalski} and Piotr~M. {So{\l}tan}.
\newblock {Quantum families of invertible maps and related problems}.
\newblock {\em Canad. J. Math.}, 68(3):698--720, 2016.

\bibitem[Tim08]{timmermann}
Thomas Timmermann.
\newblock {\em An invitation to quantum groups and duality}.
\newblock EMS Textbooks in Mathematics. European Mathematical Society (EMS),
  Z\"urich, 2008.
\newblock From Hopf algebras to multiplicative unitaries and beyond.

\bibitem[VDW96]{vDW}
Alfons Van~Daele and Shuzhou Wang.
\newblock Universal quantum groups.
\newblock {\em Internat. J. Math.}, 7(2):255--263, 1996.

\bibitem[Voi86]{Voi86}
Dan Voiculescu.
\newblock Addition of certain noncommuting random variables.
\newblock {\em J. Funct. Anal.}, 66(3):323--346, 1986.

\bibitem[Wor87]{SLW87a}
Stanis{\l}aw~Lech Woronowicz.
\newblock Compact matrix pseudogroups.
\newblock {\em Comm. Math. Phys.}, 111(4):613--665, 1987.

\bibitem[Wor98]{SLW95b}
Stanis{\l}aw~Lech Woronowicz.
\newblock Compact quantum groups.
\newblock In {\em Sym\'etries quantiques ({L}es {H}ouches, 1995)}, pages
  845--884. North-Holland, Amsterdam, 1998.

\end{thebibliography}

\end{document}